\documentclass[12pt]{article}
\usepackage{amsmath,amssymb,amsbsy,amsfonts,amsthm,latexsym,
                        amsopn,amstext,amsxtra,euscript,amscd}

\begin{document}

%%%%%%%%%%%%%%%%%%%%%%%%%%%%%%%%%%%%%%%%%%%%%%%%%%%%%%%%
%%%%%%%%%%%%%%%%%%%%%%%%%%%%%%%%%%%%%%%%%%%%%%%%%%%%%%%%
%%%%%%%%%%%%%%%%%%%%%%%%%%%%%%%%%%%%%%%%%%%%%%%%%%%%%%%%
%%%%%%%%%%%%%%%%%%%%%%%%%%%%%%%%%%%%%%%%%%%%%%%%%%%%%%%%

%%%%%%%   STANDARD STUFF %%%%%%%%%

%%%%%%%%%%%%%%%%%%%%%%%%%%%%%%%%%%%%%%%%%%%%%%%%%%%%%%%%
%%%%%%%%%%%%%%%%%%%%%%%%%%%%%%%%%%%%%%%%%%%%%%%%%%%%%%%%
%%%%%%%%%%%%%%%%%%%%%%%%%%%%%%%%%%%%%%%%%%%%%%%%%%%%%%%%
%%%%%%%%%%%%%%%%%%%%%%%%%%%%%%%%%%%%%%%%%%%%%%%%%%%%%%%%

%  use the AMS-Euler Fraktur fonts
%%%%%%%%%%%%%%%%%%%%%%%%%%%%%%%%%%
\newfont{\teneufm}{eufm10}
\newfont{\seveneufm}{eufm7}
\newfont{\fiveeufm}{eufm5}
%%%%%%%%%%%%%%%%%%%%%%%%%%%%%%%%%
%
%  allow automatic size selection in math mode
%
%%%%%%%%%%%%%%%%%%%%%%%%%%%%%%%%%
\newfam\eufmfam
      \textfont\eufmfam=\teneufm \scriptfont\eufmfam=\seveneufm
      \scriptscriptfont\eufmfam=\fiveeufm
%%%%%%%%%%%%%%%%%%%%%%%%%%%%%%%%%
%
%  \frak works on a single symbol at a time\ldots
%
\def\frak#1{{\fam\eufmfam\relax#1}}
%

%%%%%%%%%%%%%%%%%%%  bbb-matter

\def\bbbr{{\rm I\!R}} %reelle Zahlen
\def\bbbm{{\rm I\!M}}
\def\bbbn{{\rm I\!N}} %natuerliche Zahlen
\def\bbbf{{\rm I\!F}}
\def\bbbh{{\rm I\!H}}
\def\bbbk{{\rm I\!K}}
\def\bbbp{{\rm I\!P}}
\def\bbbone{{\mathchoice {\rm 1\mskip-4mu l} {\rm 1\mskip-4mu l}
{\rm 1\mskip-4.5mu l} {\rm 1\mskip-5mu l}}}
\def\bbbc{{\mathchoice {\setbox0=\hbox{$\displaystyle\rm C$}\hbox{\hbox
to0pt{\kern0.4\wd0\vrule height0.9\ht0\hss}\box0}}
{\setbox0=\hbox{$\textstyle\rm C$}\hbox{\hbox
to0pt{\kern0.4\wd0\vrule height0.9\ht0\hss}\box0}}
{\setbox0=\hbox{$\scriptstyle\rm C$}\hbox{\hbox
to0pt{\kern0.4\wd0\vrule height0.9\ht0\hss}\box0}}
{\setbox0=\hbox{$\scriptscriptstyle\rm C$}\hbox{\hbox
to0pt{\kern0.4\wd0\vrule height0.9\ht0\hss}\box0}}}}
\def\bbbq{{\mathchoice {\setbox0=\hbox{$\displaystyle\rm
Q$}\hbox{\raise
0.15\ht0\hbox to0pt{\kern0.4\wd0\vrule height0.8\ht0\hss}\box0}}
{\setbox0=\hbox{$\textstyle\rm Q$}\hbox{\raise
0.15\ht0\hbox to0pt{\kern0.4\wd0\vrule height0.8\ht0\hss}\box0}}
{\setbox0=\hbox{$\scriptstyle\rm Q$}\hbox{\raise
0.15\ht0\hbox to0pt{\kern0.4\wd0\vrule height0.7\ht0\hss}\box0}}
{\setbox0=\hbox{$\scriptscriptstyle\rm Q$}\hbox{\raise
0.15\ht0\hbox to0pt{\kern0.4\wd0\vrule height0.7\ht0\hss}\box0}}}}
\def\bbbt{{\mathchoice {\setbox0=\hbox{$\displaystyle\rm
T$}\hbox{\hbox to0pt{\kern0.3\wd0\vrule height0.9\ht0\hss}\box0}}
{\setbox0=\hbox{$\textstyle\rm T$}\hbox{\hbox
to0pt{\kern0.3\wd0\vrule height0.9\ht0\hss}\box0}}
{\setbox0=\hbox{$\scriptstyle\rm T$}\hbox{\hbox
to0pt{\kern0.3\wd0\vrule height0.9\ht0\hss}\box0}}
{\setbox0=\hbox{$\scriptscriptstyle\rm T$}\hbox{\hbox
to0pt{\kern0.3\wd0\vrule height0.9\ht0\hss}\box0}}}}
\def\bbbs{{\mathchoice
{\setbox0=\hbox{$\displaystyle     \rm S$}\hbox{\raise0.5\ht0\hbox
to0pt{\kern0.35\wd0\vrule height0.45\ht0\hss}\hbox
to0pt{\kern0.55\wd0\vrule height0.5\ht0\hss}\box0}}
{\setbox0=\hbox{$\textstyle        \rm S$}\hbox{\raise0.5\ht0\hbox
to0pt{\kern0.35\wd0\vrule height0.45\ht0\hss}\hbox
to0pt{\kern0.55\wd0\vrule height0.5\ht0\hss}\box0}}
{\setbox0=\hbox{$\scriptstyle      \rm S$}\hbox{\raise0.5\ht0\hbox
to0pt{\kern0.35\wd0\vrule height0.45\ht0\hss}\raise0.05\ht0\hbox
to0pt{\kern0.5\wd0\vrule height0.45\ht0\hss}\box0}}
{\setbox0=\hbox{$\scriptscriptstyle\rm S$}\hbox{\raise0.5\ht0\hbox
to0pt{\kern0.4\wd0\vrule height0.45\ht0\hss}\raise0.05\ht0\hbox
to0pt{\kern0.55\wd0\vrule height0.45\ht0\hss}\box0}}}}
\def\bbbz{{\mathchoice {\hbox{$\sf\textstyle Z\kern-0.4em Z$}}
{\hbox{$\sf\textstyle Z\kern-0.4em Z$}}
{\hbox{$\sf\scriptstyle Z\kern-0.3em Z$}}
{\hbox{$\sf\scriptscriptstyle Z\kern-0.2em Z$}}}}
\def\ts{\thinspace}

\newtheorem{theorem}{Theorem}
\newtheorem{lemma}[theorem]{Lemma}
\newtheorem{claim}[theorem]{Claim}
\newtheorem{cor}[theorem]{Corollary}
\newtheorem{prop}[theorem]{Proposition}
\newtheorem{definition}{Definition}
\newtheorem{question}[theorem]{Open Question}

\def\squareforqed{\hbox{\rlap{$\sqcap$}$\sqcup$}}
\def\qed{\ifmmode\squareforqed\else{\unskip\nobreak\hfil
\penalty50\hskip1em\null\nobreak\hfil\squareforqed
\parfillskip=0pt\finalhyphendemerits=0\endgraf}\fi}

%%%%%%%%%%%%%%%%%%%%%%%%%%%%%%%%%%%%%%%%%%%%%%%%%%%%%%%%
%%%%%%%%%%%%%%%%%%%%%%%%%%%%%%%%%%%%%%%%%%%%%%%%%%%%%%%%
%%%%%%%%%%%%%%%%%%%%%%%%%%%%%%%%%%%%%%%%%%%%%%%%%%%%%%%%
%%%%%%%%%%%%%%%%%%%%%%%%%%%%%%%%%%%%%%%%%%%%%%%%%%%%%%%%

%%%%%%%  END OF STANDARD STUFF %%%%%%%%%

%%%%%%%%%%%%%%%%%%%%%%%%%%%%%%%%%%%%%%%%%%%%%%%%%%%%%%%%
%%%%%%%%%%%%%%%%%%%%%%%%%%%%%%%%%%%%%%%%%%%%%%%%%%%%%%%%
%%%%%%%%%%%%%%%%%%%%%%%%%%%%%%%%%%%%%%%%%%%%%%%%%%%%%%%%
%%%%%%%%%%%%%%%%%%%%%%%%%%%%%%%%%%%%%%%%%%%%%%%%%%%%%%%

\newcommand{\ignore}[1]{}

\newcommand{\comm}[1]{\marginpar {\fbox{#1}}}

\hyphenation{re-pub-lished}

\baselineskip 15pt

\def\lln{{\mathrm Lnln}}
\def\ad{{\mathrm ad}}

\def\vec#1{\mathbf{#1}}

\def \F{{\bbbf}}
\def \K{{\bbbk}}
\def \Z{{\bbbz}}
\def \N{{\bbbn}}
\def \Q{{\bbbq}}
\def \R{{\bbbr}}
\def\Fp{\F_p}
\def \fp{\Fp^*}

\def\cA{{\mathcal A}}
\def\cE{{\mathcal E}}
\def\cI{{\mathcal I}}
\def\cQ{{\mathcal Q}}
\def\cR{{\mathcal R}}
\def\cX{{\mathcal X}}
\def\cU{{\mathcal U}}
\def\Zm{\Z_m}
\def\Zt{\Z_t}
\def\Zp{\Z_p}
\def\Fp{\F_p}
\def\Um{\Uc_m}
\def\Ut{\Uc_t}
\def\Up{\Uc_p}

\def\Ikl{\widetilde{I}_\ell(K,H,N)}
\def\\{\cr}
\def\({\left(}
\def\){\right)}
\def\fl#1{\left\lfloor#1\right\rfloor}
\def\rf#1{\left\lceil#1\right\rceil}
\def\ep{\mathbf{e}}

\newcommand{\vt}{\vartheta}
\newcommand{\floor}[1]{\lfloor {#1} \rfloor }

%%%%%%%%%%%%%%%  Topmatter %%%%%%%%%%%%%%%%%%

\title{Exponential Sums and Congruences with Factorials}

\author{
{Moubariz~Z.~Garaev}\\
\normalsize{Instituto de Matem{\'a}ticas,  Universidad Nacional Aut\'onoma de
M{\'e}xico}
\\
\normalsize{C.P. 58180, Morelia, Michoac{\'a}n, M{\'e}xico} \\
\normalsize{\tt garaev@matmor.unam.mx} \\
\and
{ Florian~Luca} \\
\normalsize {Instituto de Matem{\'a}ticas, Universidad Nacional Aut\'onoma de
M{\'e}xico}
\\
\normalsize{C.P. 58180, Morelia, Michoac{\'a}n, M{\'e}xico} \\
\normalsize{\tt fluca@matmor.unam.mx} \\
\and
{Igor E.~Shparlinski} \\
\normalsize{Department of Computing, Macquarie University} \\
\normalsize{Sydney, NSW 2109, Australia} \\
\normalsize{\tt igor@ics.mq.edu.au}}
\date{\today}

\pagenumbering{arabic}

\maketitle

%\newpage

\begin{abstract} We estimate the number of solutions of  certain
diagonal congruences involving factorials. We use these  results to bound
exponential sums
with products of two factorials $n!m!$ and also derive
asymptotic formulas for the number of solutions of various
congruences with factorials. For example, we prove that the products of
two factorials $n!m!$ with
$\max\{n,m\}<p^{1/2+\varepsilon}$ are uniformly distributed modulo
$p$, and that
any  residue class modulo $p$ is representable in the form $m!n!+n_1!
+ \ldots +n_{49}!$
with $\max  \{m,n, n_1, \ldots, n_{49}\} < p^{8775/8794+ \varepsilon}$.
\end{abstract}


\section{Introduction}

Throughout this paper, $p$ is an odd prime.
Very little
seems to be known, or even conjectured,
about the distribution of $n!$ modulo $p$.
In {\bf F11} in~\cite{RKG}, it is conjectured that
about $p/e$ of the residue classes $a\pmod p$ are missed
by the sequence $n!$. If this were so, the sequence
$n!$ modulo $p$ should assume about $(1-1/e)p$ distinct values. Some
results of  this spirit have appeared in~\cite{CVZ}.

The scarcity of heuristic results is probably due to the hardness of
computing factorials.
The best  know algorithm to compute $n!$ over $\Z$ or modulo $p$ takes about
$n^{1/2}$ arithmetic operations in the corresponding ring,
see~\cite{BCS, Cheng}.
It has been shown in~\cite{ShSm} that the complexity of computing
factorials is related to such deep conjectures of the complexity theory
as the algebraic version of  the {\bf P = NP} question, see also a nice
discussion in~\cite{Cheng}.

Sums of multiplicative characters and various
additive and multiplicative congruences with
factorials have been considered in~\cite{GaLu,GaLuSh,LuSt}.

In particular, it is has been  shown in~\cite{GaLuSh} that
for any nonprincipal character $\chi$ modulo $p$ we have
\begin{equation}
\label{eq: Bound Char Sum}
\sum_{n=L+1}^{L+N}\chi(n!) =O\( N^{3/4}p^{1/8}(\log p)^{1/4}\)
\end{equation}
and also that the number of solutions $I_\ell(L,N)$ of the congruence
\begin{equation}
\label{eq: Prod fact}
n_1!\ldots n_\ell!\equiv n_{\ell+1}!\ldots n_{2\ell}!\pmod{p}, \qquad L <
n_1,\ldots,n_{2\ell} \le L+N ,
\end{equation}
satisfies the bound
\begin{equation}
\label{eq: Bound I}
I_\ell(L,N) \ll  N^{2\ell-1+2^{-l}}
\end{equation}
(provided $0 \le L < L+N< p$).

Using~\eqref{eq: Bound Char Sum}  and~\eqref{eq: Bound I},
it has been shown in~\cite{GaLuSh} that for any fixed $\varepsilon > 0$
 the products
of three factorials $n_1! n_2! n_3!$, with
$\max \{n_1, n_2, n_3\}  = O( p^{5/6+\varepsilon})$ are
uniformly distributed modulo $p$. Here, we obtain an
upper bound for the additive analogue
of~\eqref{eq: Prod fact} and use it to   estimate double
exponential sums with products of two factorials. Namely,
for integers  $a$, $K$, $L$, $M$ and $N$ we consider
double exponential  sums
$$
W_a(K,M; L,N) =\sum_{m = K+1}^{K+M} \sum_{n = L+1}^{L+N}
\ep(a m! n!),
$$
where we define
$$
\ep(z) = \exp(2 \pi i z/p).
$$
In turn, our bound of exponential sums $W_a(K,M; L,N)$
lead us to  a substantial
improvement of the aforementioned result, showing that  the  products of two
factorials $n!m!$, with
$\max \{n, m\}  = O( p^{1/2+\varepsilon})$, are
uniformly distributed modulo $p$.
We then  combine  our new bounds and the
bounds~\eqref{eq: Bound Char Sum}
and~\eqref{eq: Bound I} to study various congruences
involving factorials.

%Our approach to estimating double sums $W_a(K,M; L,N)$
%dates back to work of Vinogradov on bounds of exponential
%sums with polynomials, and has occurred in the literature
%in various forms, for example, see Lemma~4 in~\cite{Kony} and
%the follow-up discussion. Nevertheless, it has never been applied to
%exponential sums involving factorials.

Studying single exponential sums
$$
S_a(L,N)= \sum_{n = L+1}^{L+N} \ep\(a n!\)
$$
is of great interest too.
Although we have not been able to obtain  ``individual'' bounds
for these sums, we obtain various bounds ``on average'', which also
play a crucial role in our arguments.

 Some results and techniques  of~\cite{GaLuSh} have found their applications
to studying prime divisors of $n!+f(n)$ for various functions $f$,
see~\cite{LuSh1,LuSh2}.
In particular, in~\cite{LuSh1} they have led to an improvement of
a result of Erd{\H o}s  and  Stewart~\cite{ErSt}.  We expect  that the results
of this work will also find some applications to various arithmetic
questions.

Throughout the paper, the  implied constants in symbols `$O$' and  `$\ll$'
may occasionally, where obvious, depend on   integer parameters
$k$, $\ell$, $r$ and  a small real parameter $\varepsilon > 0$,
 and are absolute otherwise
(we recall that $U \ll V$ and   $A = O(B)$ are both equivalent
to the inequality $|U| \le c V$ with some constant $c> 0$).

\bigskip

{\bf Acknowledgements.}
During the preparation of this paper,
F.~L.\ was supported in part by grants
SEP-CONACYT 37259-E and 37260-E, and
I.~S.\ was supported in part by ARC grant  DP0211459.

\section{Bounds on the Number of Solutions of  Additive  Congruences with
Factorials}

For    integers $\ell\ge 1 $, $\lambda$,    $L$ and $N$  with  $0 \le
L< L+ N < p$
we  denote by $J_\ell(L,N; \lambda)$ the number of solutions
to the congruence
$$
\sum_{i=1}^\ell  n_i! \equiv \sum_{i=\ell+1}^{2\ell} n_i! +\lambda \pmod p,
\qquad  L +1 \le n_1, \ldots, n_{2\ell} \le L+N.
$$
We also put $J_\ell(L,N) =  J_\ell(L,N,0)$.

Our treatment of $J_\ell(L,N; \lambda)$
is based on exponential sums. Accordingly,
we recall the identity
\begin{equation}
\label{eq:Ident}
\sum_{a=0}^{p-1} \ep( au) =
\left\{ \begin{array}{ll}
0,& \quad \mbox{if}\ u\not \equiv 0 \pmod p, \\
p,& \quad \mbox{if}\ u \equiv 0 \pmod p,
\end{array} \right.
\end{equation}
which we will repeatedly use, in particular to relate the  number of solutions
of various congruences and exponential sums.

\begin{theorem}
\label{thm:DiagEq}
Let    $L$  and $N$ be integers with   $0 \le L< L+ N < p$.
Then for any positive   integer  $\ell \ge 1$,
the inequality
$$
 J_\ell(L,N; \lambda) \ll N^{2\ell - 1 + 1/(\ell + 1) }
$$
holds.
\end{theorem}

\begin{proof} The identity~\eqref{eq:Ident} implies that
\begin{equation}
\label{eq:J_a and S}
J_\ell(L,N; \lambda) = \frac{1}{p} \sum_{a =0}^{p-1} |S_a(L,N)|^{2
\ell} \ep\(-a \lambda\) .
\end{equation}
In particular
$$
J_\ell(L,N; \lambda) \le J_\ell(L,N) = \frac{1}{p} \sum_{a =0}^{p-1}
|S_a(L,N)|^{2 \ell} .
$$

For any integer $k \ge 0$ we have
$$
S_a(L,N)= \sum_{n =L+1}^{L+N}
\ep\(a(n+k)!\)  + O(k).
$$
Therefore, for any integer $K \ge 0$,
\begin{eqnarray*}
S_a(L,N)&= & \frac{1}{K} \sum_{k = 0}^{K-1}\sum_{n =L+1}^{L+N}
     \ep\(a(n+k)!\)   + O(K) \\
 &= & \frac{1}{K}\sum_{n =L+1}^{L+N}
     \sum_{k = 1}^{K} \ep\( an!
\prod_{i=1}^k (n+ i)\) + O(K) \\ & = & \frac{1}{K}\sum_{n =L+1}^{L+N}
\sum_{k = 0}^{K-1}
\ep\( a n!\prod_{i=1}^k (n+ i)\)   + O(K).
\end{eqnarray*}

Using the H{\"o}lder inequality,
we derive
$$\sum_{a =0}^{p-1}|S_a(L,N)|^{2\ell}
\ll  K^{-2 \ell}
N^{2\ell -1} \sum_{k_1, \ldots, k_{2\ell}= 0}^{K-1}
\sum_{n =L+1}^{L+N}  \sum_{a =0}^{p-1} \ep\( a n! \Phi_{k_1, \ldots,
k_{2\ell}}(n)\) + p
K^{2\ell},
$$
where
$$
\Phi_{k_1, \ldots, k_{2\ell}}(X)  =
\sum_{\nu=1}^{\ell}\prod_{i=1}^{k_\nu} (n+ i) -
\sum_{\nu=\ell +1}^{2\ell} \prod_{i=1}^{k_\nu} (n+ i).
$$
The sum over $a$ vanishes, unless
\begin{equation}
\label{eq:Phi=1}
n!\Phi_{k_1, \ldots, k_{2\ell}}(n) \equiv  0 \pmod p,
\end{equation}
in which case it equals $p$.

It is easy to see that  $\Phi_{k_1, \ldots, k_{2\ell}}(X)$
is a nonconstant polynomial of degree
$O(K)$, unless $\(k_1, \ldots, k_{\ell}\)$ is a permutation of
$\(k_{\ell + 1}, \ldots, k_{2\ell}\)$, which happens for $O(K^\ell)$
choices of  $0 \le k_1,
\ldots, k_{2\ell} \le K-1$. If $\Phi_{k_1, \ldots, k_{2\ell}}(X)$
is a nonconstant polynomial,  then~\eqref{eq:Phi=1} is satisfied for
at most $K$ values of $n$,
otherwise we use the trivial bound $N$ on the number of solutions in $n$.

 Because $n! \not\equiv 0
\pmod p$ for
$0 \le L < n \le L + N < p$, the  total number of solutions
of~\eqref{eq:Phi=1}  in $L+1
\le n \le L+N$ and $0 \le k_1,
\ldots, k_{2\ell} \le K-1$ is $O(N K^{\ell}  + K^{2 \ell + 1})$.

Thus
\begin{eqnarray*}
\sum_{a =0}^{p-1}|S_a(L,N)|^{2\ell} &\ll &  K^{-2 \ell}
N^{2\ell -1} \(N K^{\ell}  +  K^{2 \ell + 1}\) p +
K^{2 \ell} p \\ &= & \(N^{2\ell}K^{-\ell}  +  N^{2\ell -1}  K
+ K^{2 \ell}  \)  p.
\end{eqnarray*}
Taking $K = \fl{ N^{1/(\ell + 1)}}$ and remarking that with this value
of $K$ the last term never dominates,  we  finish the proof.
\end{proof}

\begin{cor}
\label{le:AltSignCongr}
Let $\delta_i=\pm 1$ for each $ i = 1,\ldots, k$.  Then for any
integer $\lambda$, the
number of solutions of the congruence
$$
\sum_{i=1}^{k}\delta_in_i!\equiv \lambda \pmod{p},\qquad L+1\le
n_1,\ldots,n_k\le L+N,
$$
is $O\(N^{k-1+1/2(k_1+1)+ 1/2(k_2+1)}\)$,
where $k_1 = \fl{k/2}$, $k_2 = \fl{(k+1)/2}$.
\end{cor}

\begin{proof} Note that $k=k_1+k_2$. By the
identity~\eqref{eq:Ident}, the number $J$
of solutions of the above congruence can be expressed
via  exponential sums as
\begin{eqnarray*}
J & = &\frac{1}{p}
\ \sum_{n_1,\ldots,n_k=L+1}^{L+N}
\sum_{a=0}^{p-1}\ep\(a\(\sum_{i=1}^k\delta_i n_i! - \lambda\)\) \\
& = &
\frac{1}{p}\sum_{a=0}^{p-1} \ep(- a\lambda)
\prod_{i=1}^k S_{\delta_i a}(L,N).
\end{eqnarray*}
Since $|S_{\delta_i a}(L,N)| = |S_a(L,N)|$, we see that
$$
J\le\frac{1}{p}\sum_{a=0}^{p-1}
\left|S_a(L,N)\right|^{k}
= \frac{1}{p}\sum_{b=0}^{p-1}
\left|S_a(L,N)\right|^{k_1}
\left|S_a(L,N)\right|^{k_2}.
$$
Using the Cauchy inequality and Theorem~\ref{thm:DiagEq} (see
also~\eqref{eq:J_a
and S}), we finish the proof.
\end{proof}

Let $F_\ell(K,M; L,N)$ denote the number of solutions of the congruence
\begin{equation}
\label{eq:DoubleDiagEq}
\begin{split}
\sum_{i=1}^{\ell}n_i!m_i!\equiv
\sum_{i=\ell+1}^{2\ell}n_{i}!m_{i}!&\pmod{p},\\
K+1 \le  m_1,\ldots , m_{2\ell}\le K+M,   \quad &
L+1 \le n_1,  \ldots ,n_{2\ell}\le L+N.
\end{split}
\end{equation}

The condition   $N^{2} \ge M\ge N^{1/2}$,  requested in our
next result can be substantially relaxed.
However, because we are mainly interested in the ``diagonal case''
$M = N$ (for which this condition is always satisfied) and in order
to avoid some  technical complications, we  use this condition.

\begin{theorem}
\label{thm:DoublDiagEq}
Let  $K$, $L$, $M$ and $N$ be integers with  $0 \le K< K+ M < p$
and  $0 \le L< L+ N < p$.
For any positive integer $\ell\ge 1$,  such that $N^{2} \ge M\ge N^{1/2}$,
the following bound holds
$$
F_\ell(K,M; L,N)\ll M^{2\ell-1+1/2\ell} N^{2\ell - 1/2(\ell +1)}.
$$
\end{theorem}

\begin{proof} First of all we note that if we prove the above inequality for
$N^{2}\ge M\ge N$ then we are done. Indeed, if $N > M \ge N^{1/2}$ then we
have $M^2\ge N\ge M$, and the statement follows from the inequality
$$
N^{2\ell-1+1/2\ell} M^{2\ell - 1/2(\ell +1)} \ll M^{2\ell-1+1/2\ell} N^{2\ell
- 1/2(\ell +1)}.
$$ 

So, let $N^2\ge M\ge N$. We set $H=\fl {M^{1- 1/2 \ell} N^{1/2 (\ell +1)}}$.
We see that $N^{1-1/2\ell(\ell+1)} \ll H \ll M^{1-1/2\ell(\ell+1)}$.
Then, by the identity~\eqref{eq:Ident},
\begin{eqnarray*}
F_\ell(K,M; L,N)
&=&\frac{1}{p}\sum_{a=0}^{p-1}\left|\sum_{m=K+1}^{K+M}
\sum_{n=L+1}^{L+N}e(am!n!)\right|^{2\ell}\\
&=&
\frac{1}{p}\sum_{a=0}^{p-1}\left|\sum_{r=1}^{H}
\sum_{K+(r-1)M/H<m\le
K+rM/H}\ \sum_{n=L+1}^{L+N}e(am!n!)\right|^{2\ell}.
\end{eqnarray*}
Applying the H{\"o}lder inequality, we obtain
\begin{eqnarray*}
\lefteqn{F_\ell(K,M; L,N)}\\
&& \qquad \ll
H^{2\ell-1}\sum_{r=1}^{H}\frac{1}{p}\sum_{a=0}^{p-1}\left|
\sum_{K+(r-1)M/H<m\le K+rM/H}\
\sum_{n=L+1}^{L+N}e(am!n!)\right|^{2\ell}.
\end{eqnarray*}
Therefore,
$$
F_\ell(K,M; L,N)\ll H^{2\ell-1} Q,
$$
where $Q$ is the number of solutions of the
congruence~\eqref{eq:DoubleDiagEq} with the additional condition that
$|m_i-m_j|\le M/H$,  $1\le i  < j \le 2\ell$.
Without loss of generality, we may suppose that
$m_1=\min\{ m_i \  | \ 1 \le i \le 2 \ell\}$.
We denote $m_1=m$ and put
$$
m_i=m+s_i, \qquad  2\le i\le 2\ell.
$$

Then the congruence~\eqref{eq:DoubleDiagEq}
in the new variables takes the form
\begin{equation}
\label{eq:NewVar}
n_1!+\sum_{i=2}^{\ell}f(m, s_i)n_i! - \sum_{i=\ell+1}^{2\ell} f(m,
s_i)n_{i}!\equiv 0 \pmod{p},
\end{equation}
where  $f(m, t)=(m+1)\ldots(m+t)$ for an
integer $t\ge 1$ and $f(n, 0)=1$.

The number of solutions of the congruence~\eqref{eq:NewVar} is collected from
two sets of variables $m$ and $n_i,s_i$, $1 \le i \le 2 \ell$:
\begin{itemize}

\item[(i)]  
the first set is such that the left hand side of~\eqref{eq:NewVar} is a
polynomial of $m$ of degree greater than zero (but less than $M/H$);

\item[(ii)]
the second set consists of those for which the
left hand side of~\eqref{eq:NewVar} is constant as a polynomial of $m$.
\end{itemize}

The number of solutions $Q_1$ of~\eqref{eq:NewVar} corresponding to
the first set
is  at most
$$
Q_1 \le N^{2\ell}\left(\frac{M}{H}+1\right)^{2\ell-1}\frac{M}{H}\ll
\(MN/H\)^{2 \ell}.
$$

For the second  set of variables, we have that as a polynomial, the left hand
side of~\eqref{eq:NewVar} is a constant. Let us numerate
$s_2,\ldots,s_{2\ell}$ in an
increasing order. Then, instead of equation~\eqref{eq:NewVar} we consider the
equation
$$
n_1!+\delta_i\sum_{i=2}^{2\ell}f(m, r_i)n_i! \equiv 0 \pmod{p},
$$
with $\delta_i=\pm 1$,  and such that
$$
0 = r_1 \le  \ldots \le r_{2 \ell} \le M/H.
$$

Moreover, for each positive integer $k \le 2 \ell$ and
positive integers $e_1, \ldots, e_k$, we
consider solutions with
\begin{equation}
\label{eq:Blocks}
\begin{split}
0 =r_1 =\ldots=r_{e_1} &<  r_{e_1+1}=\ldots=r_{e_1+e_2} <
\ldots \\
&  <  r_{e_1+\ldots+e_{k-1}+1}=\ldots=r_{e_1+\ldots+e_k}\le
M/H.
\end{split}
\end{equation}
In this case, the vanishing  of the polynomial in $m$
on the left hand side of~\eqref{eq:NewVar} leads to the
conditions
\begin{equation}
\label{eq:SystCong}
\begin{split}
n_1!+\delta_2n_{2}!\ldots+\delta_{e_1} n_{e_1}!&\equiv 0\pmod{p},\\
\delta_{e_1+1}n_{e_1+1}!+\ldots
+\delta_{e_1+e_2}n_{e_1+e_2}!&\equiv 0\pmod{p},\\
\ldots  \qquad \qquad &  \qquad \qquad\ldots \\
\delta_{e_1+\ldots+e_{k-1}+1}n_{e_1+\ldots+e_{k-1}+1}!+\ldots
+\delta_{e_1+\ldots+e_k}n_{e_1+\ldots+e_k}! &\equiv
0\pmod{p}.
\end{split}
\end{equation}
Certainly, for each solution to the   system of
congruences~\eqref{eq:SystCong} 
 there are at most $M$ possible values for $m$. We also
note that  from ~\eqref{eq:SystCong} it follows that  $e_i\ge 2$ for $1\le i\le k$. 
In particular, $k \le \ell$. 

For each  $k$-dimensional vector $\vec{e} = (e_1, \ldots, e_k)$ of
positive integers such
that
$e_1+
\ldots+ e_k  = 2\ell$, there are $O\((M/H)^{k-1}\) $ possible integer vectors
$(r_1,\ldots, r_{2\ell})$ satisfying~\eqref{eq:Blocks}.
For each such fixed vector $(r_1, \ldots, r_{2\ell})$,
the number of solutions  of the system of congruences~\eqref{eq:SystCong}, by
Corollary~\ref{le:AltSignCongr}, is at most
$$
O\( \prod_{\nu=1}^k N^{r_\nu -1 +
1/2(\fl{e_i/2}+1)+1/2(\fl{(e_i+1)/2}+1)}\)= O\(N^{2\ell  - k +
\kappa(\vec{e})}\),
$$
where
$$
\kappa(\vec{e}) =
\sum_{i=1}^{k}\(\frac{1}{2(\fl{e_i/2}+1)}+\frac{1}{2(\fl{(e_i+1)/2}+1)}\) .
$$
Therefore,
$$
Q_2 \ll \max_{\vec{e}}    (M/H)^{k-1} M N^{2\ell  - k +\kappa(\vec{e})},
$$
where the maximum is taken over all    integers $1 \le k \le \ell$ and
$k$-dimensional vectors $\vec{e} = (e_1, \ldots, e_k)$ of  integers
$e_i\ge 2$,  $1\le i\le k$, with $e_1 + \ldots + e_k  = 2\ell$.

If $\ell=1$ then $k=1, \kappa(\vec{e})  = 1/2$ and
$$
Q_2 \ll
M N^{3/2}. 
$$
Therefore, in this case we have 
$$
F_\ell(K,M; L,N)  \ll M^{2} N^{2} H^{-1} +  M N^{3/2}H,
$$
and the required estimate follows from the choice of $H$.

Now, we suppose that $\ell\ge 2$.
If $k=1$, then
$$
\kappa(\vec{e})  = \frac{1}{\ell+1} .
$$
If $k\ge 2$, then trivially
$$
\kappa(\vec{e})  \le \frac{k}{2} .
$$
Hence,
$$
Q_2 \ll
M N^{2\ell - \ell/(\ell+1)} +
N^{2\ell } H\max_{2 \le k \le 2\ell }    (M/HN^{1/2})^{k} .
$$
One verifies that for our choice of $H$ and under the condition
 $N^{\ell + 1 - 1/(\ell+1)} \ge M$ (whih is always satisfied for $\ell\ge 2$
and $N^2\ge M$),  we have
$ M/HN^{1/2}\le 1$, so the term corresponding to $k=2$
dominates. Therefore,
$$
Q_2 \ll
M N^{2\ell - \ell/(\ell+1)} +
M^2 N^{2\ell  -1}  H^{-1} .
$$
Thus, putting everything together we obtain,
\begin{eqnarray*}
F_\ell(K,M; L,N) &\le & H^{2\ell-1} \( Q_1 + Q_2\) \\
& \ll &  H^{2\ell-1} \( \(MN/H\)^{2 \ell} +  M N^{2\ell - \ell/(\ell+1)} +
M^2 N^{2\ell  -1}  H^{-1} \).
\end{eqnarray*}
Since $\ell\ge 2$ then $M\le N^{3\ell/(\ell+1)}$ for  $N^2 \ge M$. Therefore
the first term always dominates the third one, and we derive
$$
F_\ell(K,M; L,N)  \ll M^{2\ell} N^{2 \ell} H^{-1} +  M N^{2\ell -
\ell/(\ell+1)}
H^{2\ell-1}.
$$
Recalling our choice of $H$, we finish the proof.
\end{proof}

\section{Bounds of Double Exponential Sums  with Factorials}

Unfortunately we are not able to estimate single sums $S_a(L,N)$,
however we obtain nontrivial bounds for double exponential sums
with factorials.
We follow some ideas of Karatsuba~\cite{Kar2, Kar3} and Korobov~\cite{Kor1,Kor2}, see
also Lemma~4 in~\cite{Kony} and
the follow-up discussion. 

\begin{theorem}
\label{thm:Double Sum}
Let  $K$, $L$, $M$ and $N$ be integers with  $0 \le K< K+ M < p$
and  $0 \le L< L+ N < p$.
Then for any  integers  $k,\ell \ge 1$,
the inequality
$$
\max_{\gcd(a, p) = 1} |W_a(K,M; L,N)| \ll  M^{1-1/2\ell(k+1)} N^{1 -
1/2k(\ell +1)}
p^{1/2kl}
$$
holds.
\end{theorem}

\begin{proof} Let  $G_\ell(L,N; \lambda)$ denote the number of solutions
to the congruence
$$
\sum_{i=1}^\ell  n_i! \equiv \lambda \pmod p,
\qquad  L +1 \le n_1, \ldots, n_{\ell} \le L+N.
$$
Clearly,
\begin{equation}
\label{eq: Av Values}
\sum_{\lambda=0}^{p-1}  G_\ell(L,N; \lambda) = N^\ell
\qquad \text{and} \qquad
\sum_{\lambda=0}^{p-1}  G_\ell(L,N; \lambda)^2 = J_{\ell}(L,N) .
\end{equation}

By the H{\"o}lder inequality, we have
\begin{eqnarray*}
\lefteqn{
|W_a(K,M; L,N)|^{\ell} \le  M^{\ell-1} \sum_{m = K+1}^{K+M} \left|\sum_{n =
L+1}^{L+N} \ep(a m! n!)\right|^{\ell} } \\
& &\  \ = M^{\ell-1} \sum_{m = K+1}^{K+M} \left| \sum_{n_1, \ldots, n_{2\ell} =
L+1}^{L+N}
\ep\(a m! (n_1! + \ldots + n_\ell!)\)\right| \\
& &\  \ = M^{\ell-1} \sum_{m = K+1}^{K+M}  \vartheta_m \sum_{n_1,
\ldots, n_{2\ell} =
L+1}^{L+N}
\ep\(a m! (n_1! + \ldots + n_\ell!)\)
\end{eqnarray*}
for some complex numbers $\vartheta_m$ with $|\vartheta_m| = 1$,
$K+1\le m \le K+M$.

Therefore,
$$
|W_a(K,M; L,N)|^{\ell}  =   M^{\ell-1} \sum_{m = K+1}^{K+M}
\sum_{\lambda=0}^{p-1}
G_\ell(L,N; \lambda) \vartheta_m\ep\(a \lambda m! \) .
$$

Applying the  H{\"o}lder inequality again, we derive
\begin{eqnarray*}
\lefteqn{
|W_a(K,M; L,N)|^{2k\ell} \le  M^{2k(\ell-1)}
\left|\sum_{\lambda=0}^{p-1} G_\ell(L,N; \lambda)^{}
\sum_{m = K+1}^{K+M} \vartheta_m\ep\(a \lambda m! \) \right|^{2k}}\\
& &\ \  \le M^{2k(\ell-1)}
\(\sum_{\lambda=0}^{p-1} G_\ell(L,N; \lambda)^{2k/(2k-1)}\) ^{2k-1}
\sum_{\lambda=0}^{p-1}\left|\sum_{m = K+1}^{K+M} \vartheta_m\ep\(a
\lambda m! \)
\right|^{2k}.
\end{eqnarray*}

Once again, by  the  H{\"o}lder inequality,
$$
\(\sum_{\lambda=0}^{p-1} G_\ell(L,N; \lambda)^{2k/(2k-1)}\) ^{2k-1}
\le \(\sum_{\lambda=0}^{p-1} G_\ell(L,N; \lambda)\)^{2k-2}
\sum_{\lambda=0}^{p-1} G_\ell(L,N; \lambda)^{2}.
$$
Using~\eqref{eq: Av Values}, we obtain
$$
\(\sum_{\lambda=0}^{p-1} G_\ell(L,N; \lambda)^{2k/(2k-1)}\) ^{2k-1}
\le N^{2\ell(k-1)}  J_{\ell}(L,N).
$$
We now have
\begin{eqnarray*}
\lefteqn{
\sum_{\lambda=0}^{p-1}
\left|\sum_{m = K+1}^{K+M} \vartheta_m \ep\(a \lambda m! \) \right|^{2k}}\\
&  & \ \ =  \sum_{\lambda=0}^{p-1}
\sum_{m_1, \ldots, m_{2k} = K+1}^{K+M}
\prod_{\nu =1}^k \vartheta_{m_\nu} \ep\(a \lambda
m_\nu! \) \prod_{\nu =k+1}^{2k} \overline{\vartheta_{m_\nu}} \ep\(-a \lambda
m_\nu! \)\\
&  & \ \ =
\sum_{m_1, \ldots, m_{2k} = K+1}^{K+M}
\prod_{\nu =1}^k \vartheta_{m_\nu} \prod_{\nu =k+1}^{2k}
\overline{\vartheta_{m_\nu}}
 \sum_{\lambda=0}^{p-1}
\ep\(a\lambda\( \sum _{\nu =1}^k   m_\nu!  -  \sum_{\nu
=k+1}^{2k}m_\nu! \)\)\\ &&
\ \  \le p  J_k(K,M).
\end{eqnarray*}
Using Theorem~\ref{thm:DiagEq}, we obtain
\begin{eqnarray*}
|W_a(K,M; L,N)|^{2k\ell} & \le  &
p M^{2k(\ell-1)} N^{2\ell(k-1)}   J_{\ell}(L,N) J_k(K,M)\\
&\ll  &p (MN)^{2k\ell} N^{-\ell/(\ell + 1)} M^{-k/(k+1)},
\end{eqnarray*}
and the desired bound follows.  
\end{proof}

For example,  for every fixed $ \varepsilon >0$, choosing sufficiently
large $k$ and $\ell$, Theorem~\ref{thm:Double Sum}
yields   a nontrivial
bound whenever $NM \ge p^{1 + \varepsilon}$.

For $K=L =0$, $N=M = p-1$, choosing $k =\ell =2$, we
obtain that the bound of Theorem~\ref{thm:Double Sum} is of the form
$O(p^{2 - 1/24})$.

It is immediate that   Theorem~\ref{thm:Double Sum}
combined with the  Erd\H os-Tur\'an relation
between the discrepancy and the appropriate exponential sums
(see~\cite{DrTi,KuipNied})
gives essentially  the same the bound (with only an extra factor $\log p$)
 on the discrepancy of the sequence of fractional parts
$$
\left\{\frac{m! n!}{p}\right\}, \qquad
K +1 \le m \le K+M, \ L +1 \le n \le L+N.
$$
As we have remarked, this improves in several directions a similar result
from~\cite{GaLuSh}.

We also remark that Theorem~\ref{thm:DoublDiagEq}
can be reformulated as an upper bound on the average value of 
sums $W_a(K,M; L,N)$ over $a =0, \ldots, p-1$. 

\section{Asymptotic Formulas for the Number of Solutions of Mixed
 Congruences  with Factorials}

Let $T_r(K,M; L,N; \lambda )$ denote the number of solutions of the congruence
\begin{equation*}
\begin{split}
\sum_{i=1}^{r}n_i!m_i!\equiv  \lambda &\pmod{p},\\
K+1 \le  m_1,\ldots , m_{r}\le K+M,   \quad &
L+1 \le n_1,  \ldots ,n_{r}\le L+N.
\end{split}
\end{equation*}

\begin{theorem}
\label{thm:DoubleWaring}
Let  $K$, $L$, $M$ and $N$ be integers with  $0 \le K< K+ M < p$
and  $0 \le L< L+ N < p$.
For any positive   integers  $k$, $\ell$, $r$ and $s$ such that  $s
\le r/2$ and
$N^2 \ge M\ge N^{1/2}$, we have
\begin{eqnarray*}
\lefteqn{ \left| T_r(K,M; L,N; \lambda ) - \frac{(MN)^{r}}{p} \right|}\\
&& \qquad \ll
 M^{r-1+1/2s- (r-2s) /2\ell(k+1) }
 N^{r- 1/2(s+1)-(r-2s) /2k(\ell+1)}p^{(r-2s)/2k\ell}.
\end{eqnarray*}
\end{theorem}

\begin{proof} Using the standard principle, we express $T_r(K,M; L,N;
\lambda )$
via exponential sums. Then
$$
T_r(K,M; L,N; \lambda) =\frac{1}{p}\sum_{a=0}^{p-1}
\(W_a(K,M; L,N)\)^{r}\ep(-a \lambda).
$$
Separating the term $(MN)^{r}/p$ corresponding to $a=0$, we obtain
\begin{eqnarray*}
\lefteqn{ \left|T_r(K,M; L,N; \lambda )-\frac{(MN)^{r}}{p}\right|\le
\frac{1}{p}\sum_{a=1}^{p-1}\left|W_a(K,M; L,N)\right|^{r}}\\
& & \qquad \ll \max_{1\le
a\le p-1}
\left|W_a(K,M; L,N)\right|^{r-2s}\frac{1}{p}
 \sum_{a=0}^{p-1}\left|W_a(K,M; L,N)\right|^{2s} \\
& & \qquad \le  F_s(K,M; L,N) \max_{1\le
a\le p-1}
\left|W_a(K,M; L,N)\right|^{r-2s}.
\end{eqnarray*}
Using Theorem~\ref{thm:DoublDiagEq} and  Theorem~\ref{thm:Double Sum}, we
finish the proof.
\end{proof}

Taking $M = N$,  $r=7$, $s=2$, $k=\ell=2$ in
Theorem~\ref{thm:DoubleWaring},
 we derive that any residue class $\lambda$ modulo $p$
has $N^{14}p^{-1}\( + O(N^{-17/12} p^{11/8})\)$ representations  in a form
$$
m_1!n_1!+m_2!n_2!+\ldots+m_7!n_7! \equiv \lambda \pmod p,
$$
with $K+1 \le m_1, n_1, \ldots, m_7, n_7 \le K+N $  (provided
$0 \le K < K+N< p$).  In particular, each $\lambda$ is represented in the
above form for $K = 0$ and some $N$ of the size $N = O\( p^{33/34}\)$.

For     integers $r \ge 0 $, $\lambda$,    $L$ and $N$  with  $0 \le
L< L+ N < p$ we  denote by $Q_r(K,M; L,N; \lambda)$
the number of
solutions of the congruence 
\begin{equation}
\label{eq:AlmostSingleWaring}
m!n!+\sum_{1\le i\le r}n_i! \equiv \lambda \pmod{p},
\end{equation}
with $K+1 \le m \le K+M$ and
$ L +1 \le n, n_1, \ldots, n_{2\ell} \le L+N$.

\begin{theorem}
\label{thm:AlmostSingleWaring}
Let  $K$, $L$, $M$ and $N$ be integers with  $0 \le K< K+ M < p$
and  $0 \le L< L+ N < p$.
Then for any  integers  $k,\ell \ge 1$,
the inequality
\begin{eqnarray*}
\lefteqn{ \left| Q_r(K,M; L,N; \lambda ) - \frac{MN^{r+1}}{p}\right| }\\
&& \qquad \ll
M^{1-1/2\ell(k+1)}
  N^{r+1/2(r_1+1)+ 1/2(r_2+1)  - 1/2k(\ell +1)} p^{1/2kl}
\end{eqnarray*}
holds, where $r_1 = \fl{r/2}$, $r_2 = \fl{(r+1)/2}$.
\end{theorem}

\begin{proof} We express $Q_r(K,M; L,N; \lambda)$ via exponential
sums:
$$
Q_r(K,M; L,N; \lambda) = \frac{1}{p}\sum_{a=0}^{p-1}
W_a(K,M; L,N) \(S_a(L,N)\)^{r}\ep(-a \lambda).
$$
Selecting the main term $MN^{r+1}/p$, corresponding to $a =0$,
we obtain
\begin{eqnarray*}
\lefteqn{ \left| Q_r(K,M; L,N; \lambda )-\frac{MN^{r+1}}{p} \right|}\\
& & \qquad \ll \max_{1\le
a\le p-1}
\left|W_a(K,M; L,N)\right|\frac{1}{p}
 \sum_{a=0}^{p-1}\left|S_a(L,N)\right|^{r}.
\end{eqnarray*}
Using Theorem~\ref{thm:Double Sum} and the same arguments as
in the proof of Corollary~\ref{le:AltSignCongr}, we  conclude the proof.
\end{proof}

Taking $M = N$,  $r=49$, $k=\ell=2$ in
Theorem~\ref{thm:AlmostSingleWaring},
 we derive that any residue class $\lambda$ modulo $p$
has $N^{51}p^{-1}\(1  + O(N^{-4397/3900} p^{9/8})\)$ representations  in 
the form
$$
m!n!+n_1!+\ldots+n_{49}! \equiv \lambda \pmod p,
$$
with $K+1 \le m, n, n_1, \ldots,  n_{49}\le K+N $  (provided
$0 \le K < K+N< p$).  In particular, each $\lambda$ is represented in the
above form for $K = 0$ and some $N$ of the size $N = O\( p^{8775/8794}\)$.

One can also derive that for any $\varepsilon> 0$ there exists an integer $r$,
such that any residue class $\lambda$ modulo $p$  can be 
representated  in the form~\eqref{eq:AlmostSingleWaring}
with $ 1\le m \le p^\varepsilon$ and  $1 \le n, n_1, \ldots, n_r<p$.

We now  combine Theorem\ref{thm:DiagEq} with the  estimate~\eqref{eq: Bound
Char Sum} from our work~\cite{GaLuSh} and apply the  method
Karatsuba~\cite{Kar1} of solving multiplicative ternary problems.

For     integers $k, \ell, r \ge 0 $,  $\lambda$,    $L$ and $N$  with  $0 \le
 L< L+ N < p$, we  denote by $R_{k,\ell,r}(K,L,S;M,N, T;\lambda)$
the number of
solutions of the congruence
$$
(m_1!+\ldots+m_k!)(n_{1}!+\ldots+n_\ell!)t_1!\ldots t_r! \equiv
\lambda \pmod{p},
$$
with  $K+1 \le  m_1,\ldots , m_{k}\le K+M$,
$L+1 \le n_1,  \ldots ,n_{\ell}\le L+N$ and
$ S+1 \le   t_1, \ldots, t_r \le S+T$.

\begin{theorem}
\label{thm:MixCong-1}
Let  $K$, $L$,   $M$, $N$, $S$, and $T$, $\lambda$ be integers with
$0 \le K< K+ M < p$,   $0 \le L< L+ N < p$ and  $0 \le S< S+ T < p$
and $\lambda \not \equiv 0 \pmod p$.
Then for any  integers  $k,\ell, r \ge 1$, the following bound holds:
\begin{eqnarray*}
\lefteqn{
R_{k,\ell,r}(K,L,S;M,N, T;\lambda)}\\
& & \qquad \qquad \ll
M^{k -1/2 + 1/2(k+1)} N^{\ell -1/2 + 1/2(\ell +1)} T^{3r/4}
p^{r/8}(\log  p)^{r/4}.
\end{eqnarray*}
\end{theorem}

\begin{proof} Let $\cX$ be the set of multiplicative characters modulo $p$,
see~\cite{LN}. We have an analogue of~\eqref{eq:Ident}
\begin{equation}
\label{eq:IdentChar}
\sum_{\chi \in \cX} \chi(u)  =
\left\{ \begin{array}{ll}
0,& \quad \mbox{if}\ u\not \equiv 1\pmod p, \\
p-1,& \quad \mbox{if}\ u \equiv 1 \pmod p.
\end{array} \right.
\end{equation}
Therefore,  we have
\begin{eqnarray*}
\lefteqn{
R_{k,\ell,r}(K,L,S;M,N, T;\lambda)}\\
& & \quad
 = \sum_{m_1,\ldots,m_k=K+1}^{K+M}
 \sum_{n_1,\ldots,n_\ell=L+1}^{L+N}  \sum_{t_1,\ldots,t_r=S+1}^{S+T}\\
& & \qquad \qquad
 \frac{1}{p-1}\sum_{\chi\in
\cX}\chi\left((m_1!+\ldots+m_k!)(n_{1}!+\ldots+n_\ell!)t_1!\ldots t_r!
\lambda^{-1}\right) \\
& & \quad
 =  \frac{1}{p-1}\sum_{\chi\in
\cX}\chi(\lambda^{-1}) \(\sum_{t =S+1}^{S+T} \chi(t)\)^r \\
& & \qquad \qquad \sum_{m_1,\ldots,m_k=K+1}^{K+M}
\chi(m_1!+\ldots+m_k!)
 \sum_{n_1,\ldots,n_\ell=L+1}^{L+N}
 \chi(n_{1}!+\ldots+n_\ell!) .
\end{eqnarray*}
Separating  the term $M^kN^\ell T^r/(p-1)$ corresponding to the
principal character
$\chi_0$ and then using the bound~\eqref {eq: Bound Char Sum}
for  the sum   over $t$,  we obtain
\begin{eqnarray*}
\lefteqn{
\left|R_{k,\ell,r}(K,L,S;M,N, T;\lambda)-\frac{M^kN^\ell T^r}{p-1}\right|  }\\
& & \quad \ll T^{3r/4}p^{r/8}(\log
 p)^{r/4}\frac{1}{p-1}\sum_{\substack{\chi \in
\cX\\\chi \ne \chi_0}}\left|\sum_{m_1,\ldots,m_k=K+1}^{K+M}
 \chi(m_1!+\ldots+m_k!)\right|\\
& & \qquad \qquad\qquad\qquad\qquad \qquad\qquad\qquad\left|
\sum_{n_1,\ldots,n_\ell=L+1}^{L+N}\chi(n_{1}!+\ldots+n_\ell!)\right|.
\end{eqnarray*}
We see, from~\eqref{eq:IdentChar}, that
$$
\frac{1}{p-1}\sum_{\chi \in \cX}\left|\sum_{m_1,\ldots,m_k=K+1}^{K+M}
 \chi(m_1!+\ldots+m_k!)\right|^ 2 \le J_{k}(K, M)
$$
(the above estimate is almost an equality
were it not for neglecting the terms with
$m_1!+\ldots+m_k! \equiv 0 \pmod p$).  Similarly,
$$
\frac{1}{p-1}\sum_{\chi \in \cX}\left|\sum_{n_1,\ldots,n_\ell=L+1}^{L+N}
 \chi(n_1!+\ldots+n_\ell!)\right|^ 2 \le  J_{\ell}(L, N).
$$
Therefore, by the Cauchy  inequality, we see that
\begin{eqnarray*}
\lefteqn{
\left|R_{k,\ell,r}(K,L,S;M,N, T;\lambda)-\frac{M^kN^\ell T^r}{p-1}\right| }\\
& & \qquad \qquad \ll
T^{3r/4}p^{r/8}(\log  p)^{r/4}
\left(J_k(K, M)  J_{\ell}(L, N)\right)^{1/2}.
\end{eqnarray*}
Using Theorem~\ref{thm:DiagEq}, we finish the proof.
\end{proof}

Taking   $M = N = T$ ,   $k = \ell = 2$ and  $r=3$  in
Theorem~\ref{thm:MixCong-1}, we
have that any $\lambda \not \equiv 0 \pmod p$
has $N^7(p-1)^{-1} \(1  + O(  N^{-17/12} p^{11/8}(\log p)^{3/4} )\)$
representations  of
the form
$$
(m_1! + m_2!)(n_1! + n_2!) t_1!t_2!t_3! \equiv \lambda \pmod p
$$
with $K+1 \le m_1,  m_2, n_1, n_2, t_1,t_2,t_3 \le K+N$
(provided  $0 \le K < K+N< p$). The above asymptotic formula is
nontrivial for any fixed $\varepsilon > 0$ and $N \ge p^{33/34 + \varepsilon}$.

One can also take $k = 1$, $\ell = 2$,  $r = 4$ in
Theorem~\ref{thm:MixCong-1},  and
obtain an asymptotic formula for the number of representations
of the form
\begin{equation}
\label{eq:7Factorials}
(n_1! + n_2!) t_1!t_2!t_3!t_4! t_5! \equiv \lambda \pmod p,
\end{equation}
with  $K+1 \le n_1, n_2, t_1,t_2,t_3,t_4, t_5 \le K+N$,  which
becomes nontrivial
for   $N \ge p^{18/19 + \varepsilon}$. However, our next result provides a
stronger bound.

Let us define  $R_{\ell,r}(L,S;N, T;\lambda)$ as
the number of
solutions of the congruence
$$
(n_{1}!+\ldots+n_\ell!)t_1!\ldots t_r! \equiv
\lambda \pmod{p}
$$
with $L+1 \le n_1,  \ldots ,n_{\ell}\le L+N$ and
$ S+1 \le   t_1, \ldots, t_r \le S+T$.
That is,  $R_{\ell,r}(L,S;N, T;\lambda) = R_{0,\ell,r}(0,L,S;1,N, T;\lambda)$.

\begin{theorem}
\label{thm:MixCong-2}Let  $K$, $L$,   $M$, $N$, $S$, and $T$,
$\lambda$ be integers with
$0 \le K< K+ M < p$,   $0 \le L< L+ N < p$ and  $0 \le S< S+ T < p$
and $\lambda \not \equiv 0 \pmod p$.
Then for any  integers  $\ell, r\ge 1$ and $0 \le s \le r$,  the
following bound
holds:
\begin{eqnarray*}
\lefteqn{ R_{\ell,r}(L,S;N, T;\lambda)  -\frac{N^\ell T^r}{p-1}} \\
& & \qquad \ll  N^{\ell -1/2 + 1/2(\ell +1)} T^{(3r+s)/4-1/2+2^{-s-1}}
p^{(r-s)/8}(\log p)^{(r-s)/4}.
\end{eqnarray*}
\end{theorem}

\begin{proof} As in the proof of Theorem~\ref{thm:MixCong-2}, we derive
\begin{eqnarray*}
\lefteqn{
R_{\ell,r}(L,S;N, T;\lambda)  }\\
& & \quad
 =  \frac{1}{p-1}\sum_{\chi\in
\cX}\chi(\lambda^{-1}) \(\sum_{t =S+1}^{S+T} \chi(t)\)^r
 \sum_{n_1,\ldots,n_\ell=L+1}^{L+N}
 \chi(n_{1}!+\ldots+n_\ell!) .
\end{eqnarray*}
Separating  the term $N^\ell T^r/(p-1)$ corresponding to the
principal character
$\chi_0$,  and also then using the bound~\eqref {eq: Bound Char Sum}, we obtain
\begin{eqnarray*}
\lefteqn{
\left|R_{\ell,r}(L,S;N, T;\lambda)   -\frac{N^\ell T^r}{p-1}\right|  }\\
& & \qquad \ll \sum_{\substack{\chi \in
\cX\\\chi \ne \chi_0}}\left|\sum_{t =S+1}^{S+T} \chi(t)  \right|^r \left|
\sum_{n_1,\ldots,n_\ell=L+1}^{L+N}\chi(n_{1}!+\ldots+n_\ell!)\right|\\
& & \qquad \ll
T^{3(r-s)/4}p^{(r-s)/8}(\log
 p)^{(r-s)/4}\sum_{\substack{\chi \in
\cX\\\chi \ne \chi_0}}\left|\sum_{t =S+1}^{S+T} \chi(t)  \right|^s\\
& & \qquad \qquad \qquad \qquad \qquad \qquad \qquad
\left| \sum_{n_1,\ldots,n_\ell=L+1}^{L+N}\chi(n_{1}!+\ldots+n_\ell!)\right|.
\end{eqnarray*}
Applying the Cauchy inequality, as in the
proof of Theorem~\ref{thm:MixCong-2}, we derive
\begin{eqnarray*}
\lefteqn{
\left|R_{\ell,r}(L,S;N, T;\lambda)   -\frac{N^\ell T^r}{p-1}\right|  }\\
& & \qquad  \ll T^{3(r-s)/4}p^{(r-s)/8}(\log
 p)^{(r-s)/4} \(I_s(S,T)    J_\ell(L,N)\)^{1/2},
\end{eqnarray*}
 and using~\eqref{eq: Bound I} together with
Theorem~\ref{thm:DiagEq},
 we finish the proof.
\end{proof}

Taking $N=T$,  $r=5$, $\ell=2$ and  $s=2$ in Theorem~\ref{thm:MixCong-2},
 we obtain that any
$\lambda \not \equiv 0 \pmod p$ has $N^7(p-1)^{-1} \( 1 +
O(N^{-35/24}p^{11/8}(\log p)^{3/4})\)$
representations of the form ~\eqref{eq:7Factorials}
with  $K+1 \le n_1, n_2, t_1,t_2,t_3,t_4, t_5 \le K+N$ (provided
$0 \le K < K+N< p$),  which becomes
nontrivial  for   $N \ge p^{33/35 + \varepsilon}$.

\section{Remarks}

As we have remarked, a more careful examination of
the function $\kappa(\vec{e})$ would lead to a substantial relaxation
of the condition  $N^2 \ge M\ge N^{1/2}$ of
Theorem~\ref{thm:DoublDiagEq}.  This however, does not
affect the most interesting ``diagonal'' case $M=N$.

Our method  can also be used, without any substantial changes,
to study  the distribution of the  products
$$
\prod_{j=1}^n f(j) \pmod p,  \qquad L+1 \le n \le L+N,
$$
where $f(j) $ is  rational function (welldefined modulo $p$ for $j
=1, \ldots, p-1$).
In particular, with $f(j) = j^{-1}$
one can estimate exponential sums and  the number of
solutions of some  congruences with   ratios of factorials $n!/m!$.

Probably the most challenging open questions are obtaining a nontrivial
upper bound on the exponential sums $S_a(L,N)$, and also obtaining an
asymptotic
formula (or at least a lower bound) on the number of solutions
of the Waring-type congruence with factorials
$$
n_{1}!+\ldots+n_\ell!  \equiv
\lambda \pmod{p},
$$
where $L+1 \le n_1,  \ldots ,n_{\ell}\le L+N$.
Even in the case of $L=0$, $N = p-1$ these questions are still
unsolved. Theorem~\ref{thm:AlmostSingleWaring} seems to be
the closest known ``approximation'' to a full solution
of the Waring problem with factorials modulo $p$.

As we have mentioned,  our ability to  obtain  any extensive computational
evidences  is very limited.
So, we dare not make any conjectures about possible answers to the
above questions.

\end{document}